\documentclass[11pt,a4paper]{article}

\usepackage{latexsym,epsf}
\usepackage[dvips]{graphicx}
\usepackage{amsthm,amssymb,amsfonts,epsfig,color,amsmath,bbm}
\usepackage{wrapfig}

\newtheorem{theorem}{Theorem}
\newtheorem{lemma}[theorem]{Lemma}
\newtheorem{corollary}[theorem]{Corollary}

\newtheorem{conjecture}[theorem]{Conjecture}
\newtheorem{definition}[theorem]{Definition}

\ifnum\catcode`\@=\active \def\affe{@@} \else \def\affe{@} \fi

\begin{document}

\title{New results on lower bounds for the number of $(\leq k)$-facets
\thanks{Part of the research on this paper was carried out while the first author
was a visiting professor at the Mathematics Department, University
of Alcal\'a, Spain.}}

\author{
      \normalsize Oswin Aichholzer
      \footnote{Research partially supported by the FWF
      (Austrian Fonds zur F\"orderung der Wissenschaftlichen Forschung)
      under grant S09205, NFN Industrial Geometry.}\\
      \small\sf Institute for Software Technology\\
      \small\sf Graz University of Technology\\[-1mm]
      \small\sf Graz, Austria\\[-1mm]
      \small\sf {\tt oaich\affe ist.tugraz.at}
\and
      \normalsize Jes\'us Garc\'{\i}a
      \footnote{Research partially supported by
                grants MCYT TIC2002-01541, and HU2007-0017.}\\
      \small\sf Escuela Universitaria de Inform\'atica\\
      \small\sf Universidad Polit\'ecnica de Madrid\\[-1mm]
      \small\sf Madrid, Spain\\[-1mm]
      \small\sf {\tt jglopez\affe eui.upm.es}
\and
      \normalsize David Orden
      \footnote{Research partially supported by grants MTM2005-08618-C02-02, S-0505/DPI/0235-02, and HU2007-0017.}\\
      \normalsize Pedro Ramos
      \footnote{Research partially supported by grants TIC2003-08933-C02-01, S-0505/DPI/0235-02, and HU2007-0017.}\\
      \small\sf Departamento de Matem\'aticas\\
      \small\sf Universidad de Alcal\'a \\[-1mm]
      \small\sf Alcal\'a de Henares, Spain\\[-1mm]
      \small\sf {\tt [david.orden|pedro.ramos]\affe uah.es}
}

\date{}
\maketitle

\begin{abstract}
In this paper we present three different results dealing with the number of $(\leq k)$-facets
of a set of points:
\begin{enumerate}
\item We give structural properties of sets in the plane that achieve the optimal lower bound
$3\binom{k+2}{2}$ of $(\leq k)$-edges for a fixed $0\leq k\leq \lfloor n/3
\rfloor -1$;
\item We give a simple construction showing that the lower bound
$3\binom{k+2}{2}+3\binom{k-\lfloor \frac{n}{3} \rfloor+2}{2}$
for the number of $(\leq k)$-edges of a planar point set appeared in~[Aichholzer et al.
New lower bounds for the number of ($\leq k$)-edges and the rectilinear crossing number of $K_n$.
{\em Disc. Comput. Geom.} 38:1 (2007), 1--14]
is optimal in the range $\lfloor n/3 \rfloor \leq k \leq \lfloor 5n/12
\rfloor -1$;
\item We show that for $k < \lfloor n/(d+1) \rfloor$
the number of $(\leq k)$-facets of a set of $n$ points in general position in~$\mathbb{R}^d$ is at
least $(d+1)\binom{k+d}{d}$, and that this bound is tight in the given range of $k$.
\end{enumerate}
\end{abstract}

\section{Introduction}

In this paper we deal with the problem of giving lower bounds to the
number of $(\leq k)$-facets of a set of points $S$: An oriented
simplex with vertices at points of $S$ is said to be a
\emph{$k$-facet} of $S$ if it has exactly $k$ points in the positive
side of its affine hull. Similarly, the simplex is said to be an
\emph{$(\leq k)$-facet} if it has at most $k$ points in the positive
side of its affine hull. If $S\subset\mathbb{R}^2$, a $k$-facet of $S$
is usually named a $k$-edge.

The number of $k$-facets of $S$ is denoted by $e_k(S)$, and
$E_k(S)=\sum_{j=0}^k e_j(S)$ is the number of $(\leq k)$-facets (the
set $S$ will be omitted when it is clear from the context).
Giving bounds on these quantities, and on the number of the companion concept
of $k$-sets, is one of the central problems in Discrete and Computational
Geometry, and has a long history that we will not try to summarize here.
Chapter~8.3 in \cite{BMP} is a complete and up to date survey of results
and open problems in the area.

Regarding lower bounds for $E_k(S)$, which is the main topic of this
paper, the problem was first studied by Edelsbrunner et al.
\cite{ehss-89} due to its connections with the complexity of higher
order Voronoi diagrams. In that paper it was stated that, in~$\mathbb{R}^2$,
\begin{equation}\label{eq:old-lb}
E_k(S)\geq 3\binom{k+2}{2}
\end{equation}
and it was given an example showing tightness for $0\leq k\leq \lfloor n/3 \rfloor -1$.
The proof used circular sequences but, unfortunately, contained
an unpluggable gap, as pointed out by Lov\'asz et
al.~\cite{LVWW}. A correct proof, also using circular sequences,
was independently found by \'Abrego and
Fern\'andez-Merchant~\cite{af-lbrcn-05} and Lov\'asz et
al.~\cite{LVWW}. In both papers
a strong connection was discovered between the number of $(\leq k)$-edges
and the number of convex quadrilaterals in a point set $S$. Specifically, if
$\square(S)$ denotes the number of convex quadrilaterals in $S$, in~\cite{LVWW} it was shown that
\begin{equation} \label{eq:lvww}
\square(S) = \sum_{k<\frac{n-2}{2}} (n-2k-3)\,E_k(S) - \frac{3}{4}\binom{n}{3}+c_n,
\end{equation}
where
$$
\quad
c_n=
\begin{cases}
\frac{1}{4}E_{\frac{n-3}{2}}(S), &\text{if $n$ is odd},\\
0, &\text{if $n$ is even}.
\end{cases}
$$

Giving lower bounds for $\square(S)$ is
in turn equivalent to determining the rectilinear crossing number of the
complete graph: if we draw $K_n$ on top of a set of points $S$, then
the number of intersections in the drawing is exactly the number of
convex quadrilaterals in $S$. The interested reader can go through
the extensive online bibliography by Vrt'o~\cite{Vrto} where the
focus is on the problem of crossing numbers of graphs.

The lower bound in Equation~\ref{eq:old-lb} was slightly improved for $k\geq
\lfloor\frac{n}{3}\rfloor$ by Balogh and Salazar~\cite{BS}, again
using circular sequences. Using different techniques, and based on the observation
that it suffices to proof the bound for sets with triangular convex
hull, we have recently shown~\cite{AGOR-06} that, in~$\mathbb{R}^2$,
\begin{equation}\label{eq:lb}
E_k(S) \geq 3\binom{k+2}{2} + \sum_{j=\lfloor\frac{n}{3}\rfloor}^k
(3j-n+3).
\end{equation}
If $n$ is divisible by $3$, this expression can be written as
$$
E_k(S) \geq 3\binom{k+2}{2} + \binom{k-\tfrac{n}{3}+2}{2}.
$$

In this paper we deal with three different problems related to lower
bounds for $E_k$: In Section~2, we study structural properties of those
sets in~$\mathbb{R}^2$ that achieve the lower bound in Equation~\ref{eq:old-lb} for
a fixed $0\leq k\leq \lfloor n/3 \rfloor -1$. The main result of this
section is that if~$E_k(S)$ is minimum for a given~$k$, then~$E_j(S)$
is also minimum for every $0 \leq j < k$. In Section~3 we
give a construction which shows tightness of the lower bound in
Equation~\ref{eq:lb} in the range $\lfloor n/3 \rfloor
\leq k \leq \lfloor 5n/12 \rfloor -1$. Finally, in Section~4 we
study the d-dimensional version of the problem and show that, for a set
of~$n$ points in general position in~$\mathbb{R}^d$,
\begin{equation}\label{eq:lb-R3}
E_k(S) \geq (d+1)\binom{k+d}{d},\mbox{ for~$0\leq k < \lfloor \frac{n}{d+1} \rfloor$},
\end{equation}
and that this bound is tight in that range. To the best of our knowledge, this is the first result
of this kind in~$\mathbb{R}^d$.

\section{Optimal sets for $(\leq k)$-edge vectors}

Given~$S\subset\mathbb{R}^2$, let us denote by~${\cal E}_k(S)$ the set of all $(\leq k)$-edges of
$S$, hence $E_k(S)$ is the cardinality of ${\cal E}_k(S)$. Throughout this section we consider
$k \leq \lfloor \frac{n}{3} \rfloor -1$. Recall that for a fixed such~$k$, $E_k(S)$
is optimal if $E_k(S)=3{k+2 \choose 2}$. Recall also that, by definition, a $j$-edge has
exactly~$j$ points of~$S$ in the positive side of its affine hull, which in this case
is the open half plane to the right of its supporting line.

We start by giving a new, simple, and self-contained proof of the bound in Equation~\ref{eq:old-lb},
using a new technique which will be useful in the rest of the section.
Although in this section they will be used in~$\mathbb{R}^2$, the following notions are presented
in~$\mathbb{R}^d$ for the sake of generality and in view of Section~\ref{sec:Rd}.

\begin{definition}[\cite{hw-ensrq-87}] \rm
Let $S$ be a set of $n$ points and $\mathcal{H}$ a family of sets in $\mathbb{R}^d$.
A subset $N\subset S$ is
called an \emph{$\epsilon$-net} of $S$ (with respect to $\mathcal{H}$) if for every
$H\in \mathcal{H}$ such that $|H\cap S| > \epsilon n$ we have that
$H\cap N \neq \varnothing$.
\end{definition}

\begin{definition} \rm
A {\em simplicial $\epsilon$-net} of $S\subset\mathbb{R}^d$ is a set of $d+1$ vertices
of the convex hull of $S$ that are an $\epsilon$-net of $S$ with respect to
closed half-spaces. A simplicial $\tfrac{1}{2}$-net will be called a
{\em simplicial half-net}.
\end{definition}

\begin{lemma}
\label{l:net} Every set $S\subset\mathbb{R}^2$ of $n$ points has
a simplicial half-net.
\end{lemma}

\begin{proof}
Let $T$ be a triangle spanned by three vertices of the convex
hull of $S$. An edge $e$ of $T$ is called \emph{good} if the
closed half plane of its supporting line which contains the third
vertex of $T$, contains at least $\tfrac{n}{2}$ points from
$S$. $T$ is called \emph{good} if it consists of three good
edges. Clearly, the vertices of a good triangle are a simplicial
half-net of $S$.

Let $T$ be an arbitrary triangle spanned by vertices of the
convex hull of $S$ and assume that $T$ is not good. Then observe
that only one edge $e$ of $T$ is not good and let $v$ be the
vertex of $T$ not incident to $e$. Choose a point $v'$ of the
convex hull of $S$ opposite to $v$ with respect to $e$. Then $e$ and
$v'$ induce a triangle $T'$ in which $e$ is a good edge. If
$T'$ is a good triangle we are done. Otherwise we iterate this
process. As the cardinalities of the subsets of vertices of $S$ considered are strictly
decreasing (the subsets being restricted by the half plane induced by $e$), the
process terminates with a good triangle.
\end{proof}

\begin{theorem} \label{t:lb1}
For every set $S$ of $n$ points and $0 \leq k < \lfloor\frac{n-2}{2}\rfloor$ we have
$E_k(S)\geq 3\binom{k+2}{2}$.
\end{theorem}

\begin{proof}
The proof goes by induction on $n$. From Lemma~\ref{l:net}, we can guarantee the existence of $T=\{a,b,c\} \subset S$,
an $\tfrac{1}{2}$-net made up with vertices of the convex hull.

Let $S' = S \setminus T$ and
consider an edge $e\in\mathcal{E}_{k-2}(S')$. We observe that $T$ cannot
be to the right of $e$: there are at least $\tfrac{n}{2}$
points on the closed half-plane to the left of $e$ and
that would contradict the definition
of $\tfrac{1}{2}$-net. Therefore, $e\in\mathcal{E}_{k}(S)$.

If we denote by $\mathcal{ET}_k(S)$ the set of $(\leq k)$-edges of~$S$ adjacent to points
in $T$, we have that
\begin{equation}\label{eq:sets}
\mathcal{E}_{k-2}(S') \cup \mathcal{ET}_k(S)  \subset \mathcal{E}_{k}(S).
\end{equation}

There are $2(k+1)$ $(\leq k)$-edges incident to each of the convex hull
vertices $a,b,c$ (which can be obtained rotating a ray based on that vertex).
We observe that at most three edges of $\mathcal{ET}_k(S)$ might be incident to two
points of~$T$ (those of the triangle~$T$) and
that the union in Equation~\ref{eq:sets} is disjoint. Therefore,
using the induction hypothesis we have
\begin{equation}\label{eq:th1}
E_k(S) \geq E_{k-2}(S')+3+6k \geq 3\binom{k}{2} +3+6k = 3\binom{k+2}{2}.
\end{equation}
\end{proof}

\begin{corollary}
\label{cor:subset} Let $S$ be a set of $n$ points,
$T=\{a,b,c\}$ a simplicial half-net of $S$
and $S' = S \setminus T$. If $E_k(S)=3\binom{k+2}{2}$
then:
\begin{itemize}
\item[(a)] $E_{k-2}(S')=3\binom{k}{2}$.
\item[(b)] A $k$-edge of $S$ is either a $(k-2)$-edge of $S'$ or is adjacent to a point in $T$.
\end{itemize}
\end{corollary}

\begin{proof}
  If $E_k(S)=3\binom{k+2}{2}$,
  both inequalities in Equation~\ref{eq:th1} are tight. Therefore \linebreak
  $E_{k-2}(S')=3 {k \choose 2}$ and Equation~\ref{eq:sets} becomes
  $\mathcal{E}_{k-2}(S') \cup \mathcal{ET}_k(S)  = \mathcal{E}_{k}(S)$ (disjoint union) which
  trivially implies part (b).
\end{proof}

\begin{theorem}
\label{thm:hull3} If $E_k(S)=3\binom{k+2}{2}$, then $S$ has a
triangular convex hull.
\end{theorem}

\begin{figure}[htb]
   \centering
   \includegraphics[scale=0.45]{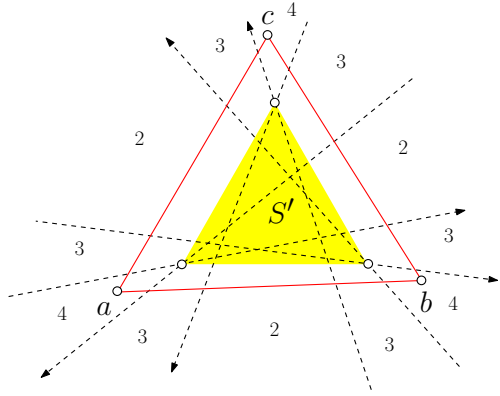}
   \caption{Each $(k-1)$-edge of $S'$ incident to a convex hull
   vertex of $S'$ (supporting lines are shown as dashed lines)
   has two vertices of $T$ on its positive side.}
   \label{fig:triangular}
\end{figure}

\begin{proof}
  We prove the statement by induction over $k$. For $k=0$ nothing has
  to be proven, so let $k=1$, assume that $E_1=9$ and let $h=|CH(S)|$.
  We have $h$ $0$-edges and at least $h$
  $1$-edges (two per convex hull vertex, but each edge might be
  counted twice). Thus $E_1 = 9 \geq 2h$ and therefore $h \leq 4$.
  Assume now $h=4$. Then at most two $1$-edges can be counted twice,
  namely the two diagonals of the convex hull. Thus we have $4+8-2=10$
  $(\leq 1)$-edges and we conclude that if $E_1=9$, then $S$ has a
  triangular convex hull.

For the general case consider $k \geq 2$, let $T=\{a,b,c\}$ be
the simplicial half-net guaranteed by Lemma~\ref{l:net}
and let $S'=S\smallsetminus T$. From Corollary~\ref{cor:subset}, part~(a),
we know that $E_{k-2}(S')=3 {k \choose 2}$ and,
by induction, we may assume that $S'$ has a triangular convex hull.
Moreover, from part~(b), no $(k-1)$-edge of $S'$ can be an $(\leq k)$-edge of $S$ and,
therefore, any $(k-1)$-edge of $S'$ must have two vertices of $T$ on its
positive side. Consider the six $(k-1)$-edges of $S'$ incident to
the three convex hull vertices of $S'$:
See Figure~\ref{fig:triangular}, where the supporting lines of these
$(k-1)$-edges are drawn as dashed lines and $S'$ is depicted as the
central triangle. Each cell outside $S'$ in the arrangement of
the supporting lines contains a number counting
the~$(k-1)$-edges considered which have that cell on their positive
side. A simple counting argument shows that the only way of placing
the three vertices $a,b,c$ of $T$ such that each $(k-1)$-edge
of $S'$ drawn has three of them on its positive side is to place one
in each cell labeled with a 4. We conclude that no vertex of~$S'$
can be on the convex hull of $S$ and the theorem follows.
\end{proof}

\begin{corollary}
\label{thm:layers3} If $E_k(S)=3\binom{k+2}{2}$, then the outermost
$\lceil \frac{k}{2} \rceil$ layers of $S$ are triangles.
\end{corollary}

\begin{proof}
{From} the optimality of $E_k(S)$ and using the same argument
as in the proof of Theorem~\ref{thm:hull3}, it follows that we
can iteratively remove the outermost $\lceil \frac{k}{2} \rceil$
layers to obtain optimal subsets, which, by Theorem~\ref{thm:hull3},
have triangular convex hulls.
\end{proof}

\begin{theorem}
\label{thm:} If $E_k(S)=3\binom{k+2}{2}$,
then $E_j(S)=3\binom{j+2}{2}$ for every $0 \leq j \leq k$.
\end{theorem}

\begin{proof}
  We prove the theorem by induction on $k$. For $k=0,1$ the theorem is
  equivalent to Theorem~\ref{thm:hull3}, so let $k \geq 2$. It is
  sufficient to show that optimality of $E_k(S)$ implies optimality of
  $E_{k-1}(S)$, as the theorem follows by induction.

  Let $T$ be the vertices of $CH(S)$ (which is a triangle as guaranteed
  by Theorem~\ref{thm:hull3}) and let $S'=S\smallsetminus T$.
  As in Theorem~\ref{t:lb1} we have
  $$
  \mathcal{E}_{k-3}(S') \cup \mathcal{ET}_{k-1}(S) \subset \mathcal{E}_{k-1}(S).
  $$
  Observe that $E_{k-2}(S')$ is optimal, as guaranteed by Corollary~\ref{cor:subset} and
  this implies optimality of $E_{k-3}(S')$  by induction. $|\mathcal{ET}_{k-1}(S)|$ is also
  optimal because the convex hull of $S$ is the triangle~$T$. Therefore, to prove optimality
  of $E_{k-1}(S)$ it only remains to show that no $(k-2)$-edge of $S'$ can be a
  $(k-1)$-edge of $S$.

  So let $e$ be a $(k-2)$-edge of $S'$ and let $p$ and $q$ be the
  vertices of the convex hull of~$S'$ incident to $e$ or on its
  positive side. The existence of $p$ and $q$ is guaranteed by
  Corollary~\ref{cor:subset}, part~(b). Without loss of generality, assume that the edge $pq$ is
  horizontal with the remaining vertices of $S'$ above it, see
  Figure~\ref{fig:bounded} for the rest of the proof. Let $\ell_1$ be
  the $(k-1)$-edge of $S'$ incident to $p$ which has $q$ on its
  positive side and $\ell_2$ the $(k-1)$-edge incident to $q$
  and having $p$ on its positive side. The \emph{boundary chain} is the
  lower envelope of~$\ell_1$,~$pq$, and $\ell_2$.
  We claim that $e$ does not intersect the boundary chain
  and lies above it. If $e$ is incident to $p$ or $q$ then the claim
  is obviously true. Otherwise observe that $e$ has to intersect the
  supporting lines of both considered $(k-1)$-edges in the interior of
  $S'$, as otherwise there would be too many vertices on the positive
  side of $e$. But then again $e$ lies above the boundary chain and
  the claim follows.

\begin{figure}[htb]
  \centering \includegraphics[scale=0.65]{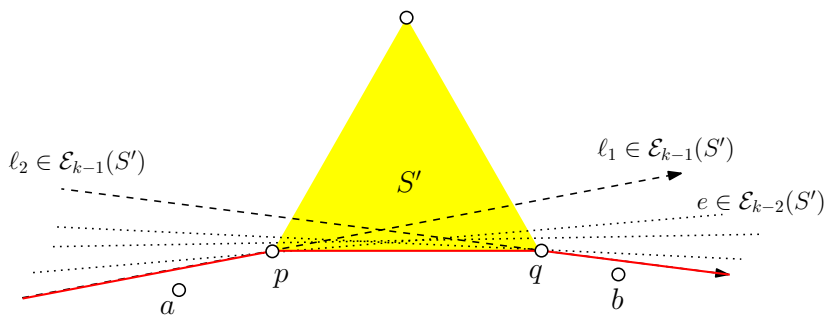}
\caption{All $(k-2)$-edges of $S'$ (supporting lines are shown as
dotted lines) lie above the (bold) lower envelope.}
\label{fig:bounded}
\end{figure}

  {From} the proof of Theorem~\ref{thm:hull3} we know that two of the
  vertices of the convex hull of $S$ have to lie below our boundary
  chain (below the $(k-1)$-edges, see $a$ and $b$ in
  Figure~\ref{fig:bounded}) and thus on the positive side of $e$.
  Therefore $e$ has at least $k$ vertices of $S$ on its positive side
  and does not belong to ${\cal E}_{k-1}(S)$. We conclude that $E_{k-1}(S)$
  is optimal and the theorem follows.
\end{proof}

\begin{corollary}\label{cor:e_j-opt}
Let $0\leq k\leq \lfloor \frac{n}{3} \rfloor-1$.
If $E_k(S)=3\binom{k+2}{2}$, then $e_j(S)=3(j+1)$ for~$0 \leq j \leq k$.
\end{corollary}


\section{Tightness of the lower bound for~$(\leq k)$-edges in~$\mathbb{R}^2$}

In this section we show a point configuration which proves tightness in the range~$0\leq
k \leq\lfloor \frac{5n}{12} \rfloor -1$ of the lower bound for~$E_k(S)$ given in~\cite{AGOR-06}.
Consider the configuration in
Figure~\ref{fig:tightness-Range1-1} (left), which is composed of
three rotationally symmetric chains, each one associated to a convex
hull vertex, fulfilling the following properties (where left and
right are considered with respect to the corresponding convex hull
vertex):

\begin{figure}[htb]
\centering
\includegraphics[width=0.9\textwidth]{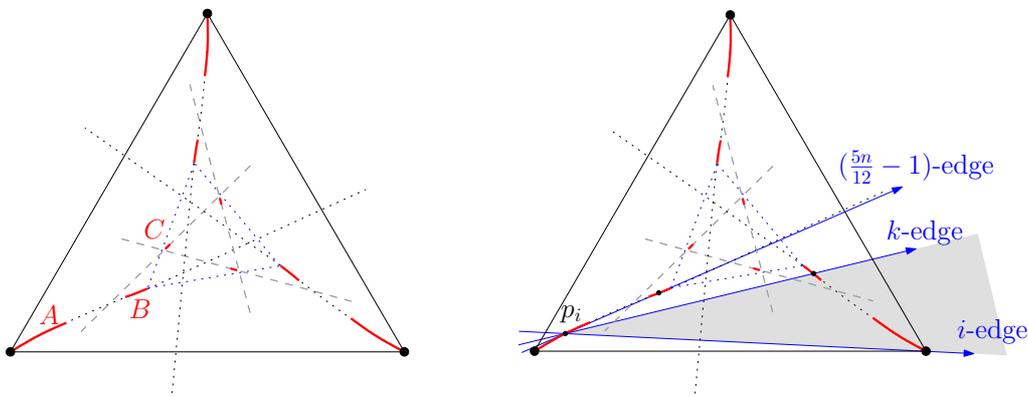} 
\caption{Left: Configuration showing tightness for~$0\leq k \leq \lfloor\frac{5n}{12}\rfloor-1$. Right:
For~$i\in\{0,\ldots,\frac{2n}{12}-1\}$, exactly one~$j$-edge appears
for each~$j\in\{i,\ldots,k\}$.} \label{fig:tightness-Range1-1}
\end{figure}

\begin{itemize}
\item The first part of the chain is slightly convex to the right and
contains~$\frac{3n}{12}$ points, with a hole between the
first~$\frac{2n}{12}$ points (which we call subchain~$A$) and the
remaining~$\frac{n}{12}$ points (called subchain~$B$).
\item Each chain is completed with a subchain~$C$, composed of
another~$\frac{n}{12}$ points slightly convex to the left.
\item All the lines spanned by two points in~$A\cup B$,
oriented from the outermost to the innermost endpoint, leave to the
right the next chain in counterclockwise order, and to the left both
the points in~$C$ and those in the remaining chain.
\item All the lines spanned by two points in~$C$ separate
subchains~$A$ and~$B$. Furthermore, when oriented from the outermost
to the innermost endpoint, they leave to the right the two remaining
subchains of type~$C$.
\item The triangle defined by the innermost points of chains of
type~$B$ contains all the chains of type~$C$.
\end{itemize}

\begin{theorem}
\label{thm:tightness} For the point configuration $S$ defined above and
$\lfloor \frac{n}{3} \rfloor\leq k \leq\lfloor \frac{5n}{12} \rfloor -1$,
\[
E_k(S)=3{{k+2}\choose {2}}+3{{k-\frac{n}{3}+2}\choose {2}},
\]
which matches the lower bound stated in~\cite{AGOR-06} for~$n$ divisible by~$3$.
\end{theorem}

\begin{proof}
Because of the rotational symmetry, we can focus on one of the three
chains~$A\cup B\cup C$ and let~$p_i$ be the $(i+1)$-th point on that
chain. We will count oriented~$j$-edges of
type~$\stackrel{\longrightarrow}{p_iq}$ (i.e. with~$p_i$ on the
tail) for~$j\leq k$. In order to do so we rotate counterclockwise a
ray based on~$p_i$, starting from the one passing through the convex
hull vertex of the next chain in counterclockwise order. Three cases
arise, depending on the index~$i$ of~$p_i$, which correspond
to~$p_i$ lying on one of the three subchains:

\begin{itemize}
\item[$(A)$] For~$i\in\{0,\ldots,\frac{2n}{12}-1\}$, exactly
one~$j$-edge appears for each~$j$ in the range~$j\in\{i,\ldots,k\}$,
while all the remaining~$j$-edges~$\stackrel{\longrightarrow}{p_iq}$
in the rotation have~$j>k$ since at some point we find
a~$(\frac{5n}{12}-1)$-edge and after that all the~$j$-edges found
have~$j\geq\frac{5n}{12}>k$. See Figure~\ref{fig:tightness-Range1-1}
(right).
\item[$(B)$] For~$i\in\{\frac{2n}{12},\ldots,\frac{3n}{12}-1\}$,
exactly one~$j$-edge appears for each~$j$ in the
ranges~$j\in\{i,\ldots,k\}$
and~$j\in\{\frac{7n}{12}-i-1,\ldots,k\}$, while all the
remaining~$j$-edges~$\stackrel{\longrightarrow}{p_iq}$ in the
rotation have~$j>k$. See Figure~\ref{fig:Range2-1and2}.
\item[$(C)$] For~$i\in\{\frac{3n}{12},\ldots,\frac{4n}{12}-1\}$,
exactly one~$j$-edge appears for each~$j$ in the
ranges~$j\in\{i,\ldots,k\}$
and~$j\in\{\frac{8n}{12}-i-1,\ldots,k\}$, while all the
remaining~$j$-edges~$\stackrel{\longrightarrow}{p_iq}$ in the
rotation have~$j>k$. See Figure~\ref{fig:Range3-1and2}.
\end{itemize}

\begin{figure}[htb]
\centering
\includegraphics[width=0.9\textwidth]{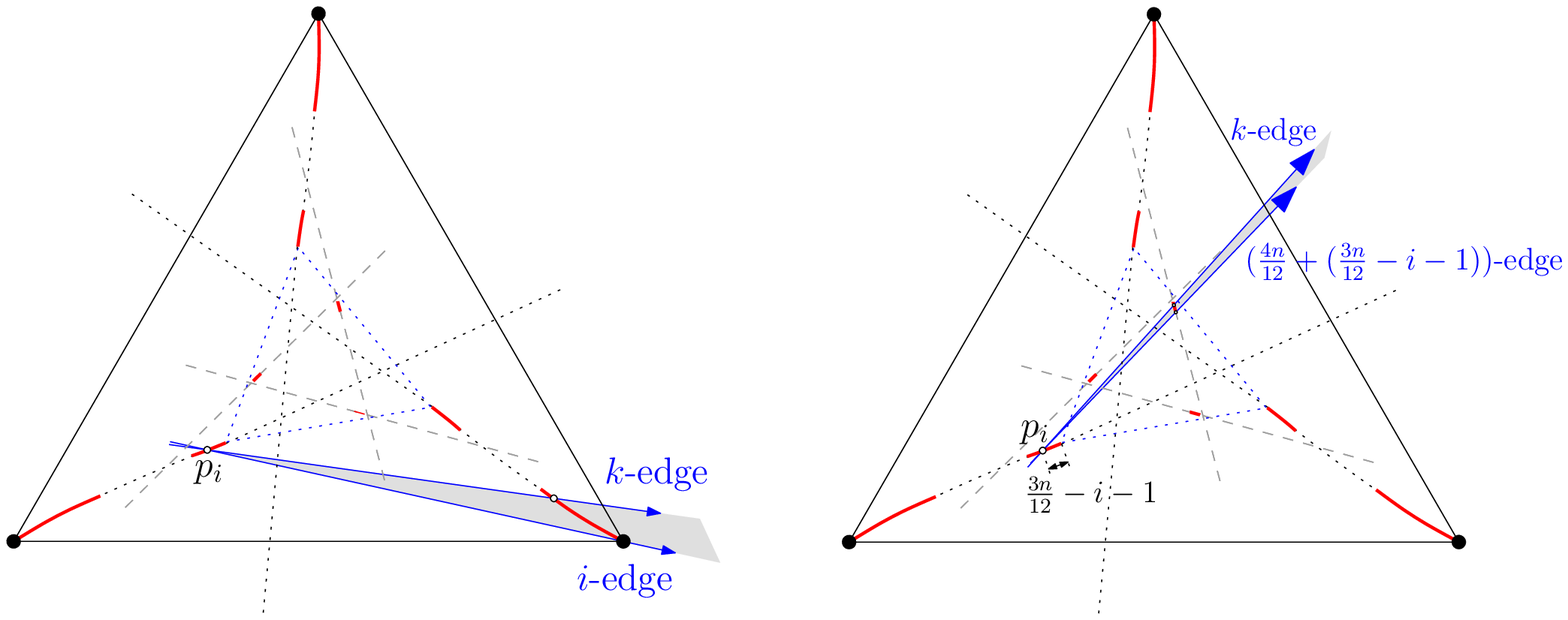} 
\caption{Left: For~$i\in\{\frac{2n}{12},\ldots,\frac{3n}{12}-1\}$,
exactly one~$j$-edge appears for each~$j\in\{i,\ldots,k\}$. Right:
For~$i\in\{\frac{2n}{12},\ldots,\frac{3n}{12}-1\}$, exactly
one~$j$-edge appears for each~$j\in\{\frac{7n}{12}-i-1,\ldots,k\}$.}
\label{fig:Range2-1and2}
\end{figure}

\begin{figure}[htb]
\centering
\includegraphics[width=0.9\textwidth]{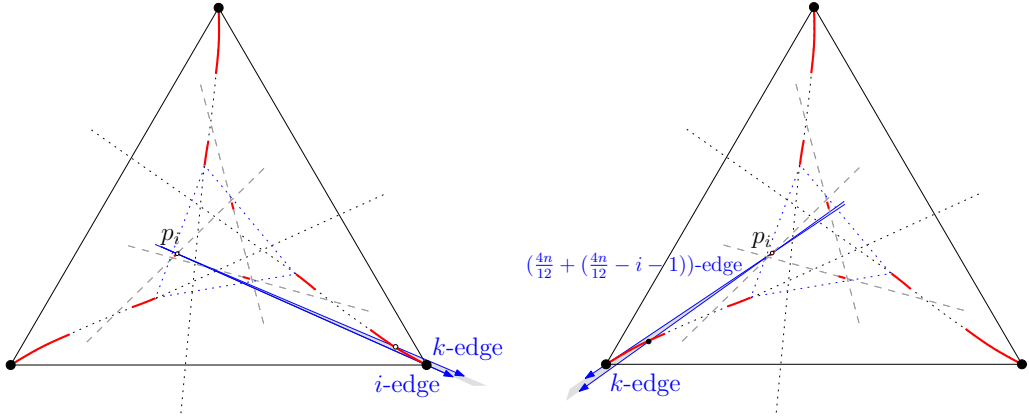} 
\caption{Left: For~$i\in\{\frac{3n}{12},\ldots,\frac{4n}{12}-1\}$,
exactly one~$j$-edge appears for each~$j\in\{i,\ldots,k\}$.
Right: For~$i\in\{\frac{3n}{12},\ldots,\frac{4n}{12}-1\}$, exactly
one~$j$-edge appears for each~$j\in\{\frac{8n}{12}-i-1,\ldots,k\}$.}
\label{fig:Range3-1and2}
\end{figure}

Let us point out that, depending on the values of~$k$ and~$i$, some
of the above ranges could actually be empty. Now we are ready to
count the total number of~$(\leq k)$-edges incident to points~$p_i$
on the first chain, which is:
\[
\sum_{i=0}^{\frac{2n}{12}-1} (k-i+1) +
\sum_{i=\frac{2n}{12}}^{\frac{3n}{12}-1} (k-i+1) +
\sum_{i=\frac{3n}{12}}^{\frac{4n}{12}-1} (k-i+1) +
\sum_{i=\frac{2n}{12}}^{\frac{3n}{12}-1} (k-\frac{7n}{12}+i+2) +
\sum_{i=\frac{3n}{12}}^{\frac{4n}{12}-1} (k-\frac{8n}{12}+i+2),
\]
where the first three summands come from the first ranges of the
three cases above, while the two remaining summands come from the
second ranges in cases~$(B)$ and~$(C)$. Merging the first three
summands and rewriting the two latter ones, the above sum equals
\[
\sum_{i=0}^{\frac{4n}{12}-1} (k-i+1) +
\sum_{i=\frac{5n}{12}-1}^{\frac{4n}{12}} (k-i+1) +
\sum_{i=\frac{5n}{12}-1}^{\frac{4n}{12}} (k-i+1) =\sum_{j=1}^{k+1} j
+ \sum_{j=1}^{k-\frac{4n}{12}+1} j = {{k+2}\choose
{2}}+{{k-\frac{n}{3}+2}\choose {2}},
\]
where the first equality comes from neglecting the negative
summands, due to the above mentioned empty ranges, and merging the
first two sums. This result has to be multiplied by the three chains
of the configuration, so we get
\[
E_k(S)=3{{k+2}\choose {2}}+3{{k-\frac{n}{3}+2}\choose {2}}.
\]
\end{proof}


\newpage

\section{A lower bound for $(\leq k)$-facets in $\mathbb{R}^d$}
\label{sec:Rd}

Throughout this section, $S\subset\mathbb{R}^d$ will be a set of~$n$
points in general position.

We remind that~$e_k(S)$ and~$E_k(S)$ denote, respectively, the
number of $k$-facets and the number of $(\leq k)$-facets of~$S$. The
main result of this section is a lower bound for the number of
$(\leq k)$-facets of a set of~$n$ points in general position in~$\mathbb{R}^d$ in the range
$0\leq k < \lfloor \tfrac{n}{d+1}\rfloor$.

The proof follows the approach in Theorem~\ref{t:lb1}, using the fact that every set of points has a centerpoint:
a point~$c\in \mathbb{R}^d$ is a \emph{centerpoint} of~$S$ if no open halfspace
that avoids~$c$ contains more than $\lceil\tfrac{dn}{d+1}\rceil$ points of~$S$ (see~\cite{Ed}).

\begin{theorem}\label{t:lbd}
Let~$S$ be a set of~$n\geq d+1$ points in~$\mathbb{R}^d$ in general
position. Then
$$
E_k(S) \geq (d+1)\binom{k+d}{d} \qquad \text{if $0\leq k < \lfloor\tfrac{n}{d+1}\rfloor$.}
$$
Furthermore, the bound on $E_k(S)$ is tight in the given range of~$k$.
\end{theorem}

\begin{proof}
The proof uses induction on~$n$ and $d$. The base case for $n=d+1$ is obvious and
for $d=2$ is just Equation~\ref{eq:old-lb}.

Let $k < \lfloor\tfrac{n}{d+1}\rfloor$ and let $c$ be a centerpoint of $S$.
Let us consider a simplex $T$ with vertices in the convex hull of $S$ and containing $c$ and
let $S'=S\smallsetminus T$.
From the definition of centerpoint, it follows that no open halfspace that avoids~$T$
contains more than \mbox{$\lceil\tfrac{dn}{d+1}\rceil-1$}~points
or, equivalently, every closed halfspace containing~$T$ has at least $\lfloor\tfrac{n}{d+1}\rfloor+1$ points.

We denote by~$\mathcal{E}_k^j(S)$ the set of
$(\leq k)$-facets of~$S$ adjacent to exactly~$j$ vertices of~$T$, and
$E_k^j(S)$ will be the cardinality of $\mathcal{E}_k^j(S)$.

For $j=0$, we observe that $\mathcal{E}_{k-d}^0(S') \subset \mathcal{E}_k^0(S)$,
because a closed halfspace containing at most $k$ points cannot contain all the
vertices of $T$. Because $k-d \leq \lfloor\tfrac{n-(d+1)}{d+1}\rfloor-1$, we can
apply induction on $n$ and get
$$
E_k^0(S)\geq E_{k-d}^0(S') \geq (d+1) \binom{(k-d)+d}{d} = (d+1) \binom{k}{d}.
$$

For $1\leq j \leq d$, let $T_j$ be a subset of $j$ vertices of $T$ and let $S_{\pi}$ be the projection
from~$T_j$ of $S\smallsetminus T$ onto the $(d-j)$-dimensional subspace $\pi$ defined by the points in
$T\smallsetminus T_j$: a point $p\in S\smallsetminus T$ is mapped to the intersection between the
$j$-flat defined by $p$ and $T_j$ and the $(d-j)$-flat defined by points in $T\smallsetminus T_j$.
Using the general position assumption, it is easy to see that the intersection has dimension zero.
If the intersection were empty, we
could slightly perturb $p$ without changing the number of $(\leq k)$-facets of $S$.

Now, if $\sigma\subset S_{\pi}$ is an $(\leq (k-d+j))$-facet of $S_{\pi}$,
then $\sigma \cup T_j$ is an $(\leq k)$-facet of $S$ (as before, a halfspace containing at most $k$ points
of $S$ cannot contain all the vertices of $T$). Because
$$
k-d+j \leq \Bigl\lfloor\frac{n}{d+1}\Bigr\rfloor -1 \leq \Bigl\lfloor\frac{n-j}{d-j+1}\Bigr\rfloor -1
$$
we can apply induction in $d$ and $n$, obtaining that there are at least
$$
(d-j+1)\binom{k-d+j+(d-j)}{d-j} = (d-j+1)\binom{k}{d-j}
$$
$(\leq k)$-facets of $S$ adjacent to $T_j$. Summing on all the subsets of $j$ points of $T$, we get
$$
E_k^j(S) \geq \binom{d+1}{j} (d-j+1)\binom{k}{d-j},
$$
and, finally,
$$
E_k(S) \geq \sum_{j=0}^d \binom{d+1}{j} (d-j+1)\binom{k}{d-j} = (d+1)\binom{k+d}{d}.
$$

As for tightness, the example showing that the bound
$3\binom{k+2}{2}$ is tight for $0\leq k\leq \lfloor\frac{n}{3}\rfloor -1$
in the planar case~\cite{ehss-89} can be extended to~$\mathbb{R}^d$: Consider $d+1$
rays in~$\mathbb{R}^d$ emanating from the origin and with the
property that any hyperplane containing one of them leaves on each open
halfspace at least one of the remaining rays. For instance, we could take the rays
defined by the origin and the vertices of a regular simplex inscribed in the unit
$d$-sphere.

Let~$n=(d+1)m$ and put chains~$C_1,\ldots,C_{d+1}$ with~$m$ points on each ray, slightly perturbed
to achieve general position. For~$j<m$, every $j$-facet of~$S$ is defined by $d$ points on different chains,
because a facet defined by two points in the same chain has at least $m$ points on each halfspace.
If we label the points of each chain from~$0$ to~$m-1$ (starting from the convex hull) and consider
$p_{i_1}^1\in C_1, \ldots,p_{i_d}^d \in C_d$, they define a $(i_1+\ldots+i_d)$-facet.
Therefore, the number of $(\leq k)$-facets defined by one point on each of these chains equals the
cardinality of the set
$$
\{(i_1,\ldots,i_d)\in\mathbb{Z}^d:i_1+\ldots+i_d\leq k, \, 0\leq i_1,\ldots,i_d\leq k\},
$$
which is exactly~$\binom{k+d}{d}$. Since these are the facets defined by points in $d$ out
of the $d+1$ chains, the total number of $(\leq k)$-facets of the set is exactly~$(d+1)\binom{k+d}{d}$.
\end{proof}


\section{Conclusions and open problems}

For~$S\subset\mathbb{R}^2$ we have shown that, for a fixed $k \leq
\lfloor \frac{n}{3} \rfloor-1$, if $E_k(S)$ is optimal, i.e. $E_k(S) =
3\binom{k+2}{2}$, then~$E_j(S)$ is also optimal in the whole range $0
\leq j \leq k$, which in turn implies that $e_j(S)=3(j+1)$ for~$0 \leq
j \leq k$. Moreover, then the outermost $\lceil \frac{k}{2}\rceil$
layers of~$S$ are triangles and these layers consist entirely of
$j$-edges of special types. In addition, we have been able to give a
simple construction showing that the lower bound in Equation~\ref{eq:lb} is
tight for $0\leq k \leq \lfloor \tfrac{5n}{12} \rfloor -1$.

All these results reveal significant deeper insight into the structure
of sets minimizing the number of $k$-edges, the final goal being to
find tight bounds for every~$k$.

Moreover, for an $n$-point set $S\subset\mathbb{R}^d$ we have proven
the lower bound $(d+1)\binom{k+d}{d}$ for the number of $(\leq k)$-facets
in the range $0\leq k < \lfloor n/(d+1) \rfloor$, which is the first
result of this kind in~$\mathbb{R}^d$.

The restriction $k < \lfloor n/(d+1) \rfloor$ stems from the underlying technique,
namely using the centerpoint of a set, and can
probably be removed. An alternative proof of Theorem~\ref{t:lbd},
using a simplicial half-net instead of a centerpoint, would be
sufficient to extend the bound to the whole range of $k$. Therefore,
it is a challenging task to extend Lemma~\ref{l:net} to dimension~$d$,
as the following conjecture states:

\begin{conjecture}\label{conj:shalfnet-Rd}
Every point set~$S\subset\mathbb{R}^d$ has a simplicial half-net.
\end{conjecture}

\end{document}